\documentclass[12pt]{amsart}
\usepackage{amssymb,amscd}
\usepackage[mathscr]{eucal}
\input epsf
\newtheorem{thm}{Theorem}[section]

\newtheorem{defin}[thm]{Definition}

\newtheorem{prop}[thm]{Proposition}
\newtheorem{rem}[thm]{Remark}

\begin{document}

\newcommand{\C}{{\Bbb C}}
\newcommand{\HH}{{\Bbb H}}
\newcommand{\M}{{\mathcal M}}
\newcommand{\R}{{\Bbb R}}
\newcommand{\Z}{{\Bbb Z}}
\newcommand{\E}{{\mathcal E}}
\newcommand{\SU}{\o{SU}(2)}
\newcommand{\Hom}{\o{Hom}}
\newcommand{\Spin}{\o{Spin}(2)}
\newcommand{\Spinm}{\o{Spin}_-(2)}
\newcommand{\Pin}{\o{Pin}(2)}
\renewcommand{\O}{{\mathcal O}}
\renewcommand{\P}{{\Bbb P}}
\newcommand{\Q}{{\mathcal Q}}
\renewcommand{\o}{\operatorname}
\newcommand{\tr}{\o{tr}}
\title
{Exceptional Discrete Mapping Class Group Orbits in Moduli Spaces}
\author{Joseph P. Previte \and Eugene Z. Xia}
\address{School of Science,
Penn State Erie, The Behrend College,
Erie, PA 16563}
\address{
National Center for Theoretical Sciences
Third General Building, National Tsing Hua University
No. 101, Sec 2, Kuang Fu Road,
Hsinchu, Taiwan 30043, Taiwan R.O.C.
}
\email{jpp@vortex.bd.psu.edu {\it (Previte)}, xia@math.umass.edu {\it (Xia)}}
\date{\today}
\maketitle
\begin{abstract}
Let $M$ be a four-holed sphere and $\Gamma$ the mapping class
group of $M$ fixing $\partial M$. The group $\Gamma$ acts on the
space  $\M_{\mathcal B}(\SU)$ of $\SU$-gauge equivalence classes
of flat $\SU$-connections on $M$ with fixed holonomy on $\partial
M$. We give examples of flat $\SU$-connections whose holonomy
groups are dense in $\SU$, but whose $\Gamma$-orbits are discrete
in $\M_{\mathcal B}(\SU)$. This phenomenon does not occur for
surfaces with genus greater than zero.
\end{abstract}

\section{Introduction}
Let $M$ be a Riemann surface of genus $g$ with $n$
boundary components (circles).  Let
$$
\{\gamma_1, \gamma_2,..., \gamma_n\} \subset \pi_1(M)
$$
be the elements in the fundamental group corresponding to
these $n$ boundary components.
Assign each $\gamma_i$ a conjugacy class $B_i \subset \SU$ and let
$$
{\mathcal B} = \{ B_1, B_2, ..., B_n \},
$$
$$
{\mathcal H}_{\mathcal B} =
\{\rho \in \Hom(\pi_1(M),\SU) : \rho(\gamma_i) \in B_i, 1 \le i \le n \}.
$$
A conjugacy class in $\SU$ is determined by its trace which is in $[-2,2]$.
Hence we might consider ${\mathcal B}$ as an element in $[-2,2]^n$.
The group $\SU$ acts on ${\mathcal H}_{\mathcal B}$
by conjugation.
\begin{defin}~\label{def:1.1}
The moduli space with fixed holonomy ${\mathcal B}$ is
$$
\M_{\mathcal B} =  {\mathcal H}_{\mathcal B}/\SU.
$$
\end{defin}
Denote by $[\rho]$ the image of $\rho \in {\mathcal H}_{\mathcal B}$
in $\M_{\mathcal B}$.
The set of smooth points of $\M_{\mathcal B}$ possesses a natural
symplectic structure which gives rise
to a finite measure $\mu$ on $\M_{\mathcal B}$ (see \cite{Go1, Go2}).

Let $\o{Diff}(M, \partial M)$ be the group of diffeomorphisms of $M$ fixing
$\partial M$.  The mapping class group $\Gamma$ is
$\pi_0(\o{Diff}(M, \partial M))$.  The group $\Gamma$
acts on $\pi_1(M)$ fixing the $B_i$'s.  This action induces a
$\Gamma$-action on $\M_{\mathcal B}$.
\begin{thm} [Goldman] \label{thm:1}
The mapping class group $\Gamma$ acts
ergodically on $\M_{\mathcal B}$ with respect to
the measure $\mu$.
\end{thm}
Since $\M_{\mathcal B}$ has a natural topology, one may also study the
topological dynamics of the mapping class group action and we have
\cite{Pr1,Pr2}:
\begin{thm} \label{thm:3}
Suppose $M$ is an orientable surface with boundary
and $g>0$.  Let
$\rho \in {\mathcal H}_{\mathcal B}$ such that $\rho(\pi_1(M))$
is dense in $\SU$.
Then the $\Gamma$-orbit of the conjugacy class
$[\rho] \in \M_{\mathcal B}$
is dense in $\M_{\mathcal B}$.
\end{thm}

In this paper we show:
\begin{thm} \label{thm:4}
Let $M$ be a four-holed sphere.  Then there exists a subset $F
\subset [-2,2]^4$ of two real dimensions with the following
property: Suppose ${\mathcal B} \in F$.  Then there exists $\rho
\in {\mathcal H}_{\mathcal B}$ with $\rho(\pi_1(M))$ dense in
$\SU$, but the $\Gamma$-orbit of the conjugacy class $[\rho]$ is
discrete in $\M_{\mathcal B}$.
\end{thm}

Let $G$ be a subgroup of $\SU$. We say that a representation
$\rho$ is a $G$-representation if $\rho(\pi_1(M)) \subset G$ up to
conjugation by $\SU$. The group $\SU$ is a double cover of
$\o{SO}(3)$:
$$
p : \SU \longrightarrow \o{SO}(3).
$$
The group $\o{SO}(3)$ contains $\o{O}(2)$,
and the symmetry groups of the regular polyhedra: ${\mathcal T}'$ (the tetrahedron),
${\mathcal C}'$ (the cube), and ${\mathcal D}'$ (the dodecahedron).
Let $\Pin, {\mathcal T}, {\mathcal C},$ and ${\mathcal D}$ denote the groups
$p^{-1}(\o{O}(2)),$ $ p^{-1}({\mathcal T}'),$ $p^{-1}({\mathcal C}'),$
and $p^{-1}({\mathcal D}')$, respectively.
The proper closed subgroups of $\SU$ consist of ${\mathcal T}$, ${\mathcal C}$, ${\mathcal D},$ and the
closed subgroups of $\Pin$.
The group $\Pin$ has two components, and
we write
$$
\Pin = \Spin \cup \Spinm,
$$
where $\Spin$ is the identity component of $\Pin$.

\begin{rem} \label{rem:dense}
Suppose $\rho \in \Hom(\pi_1(M),\SU)$.  If $\rho(\pi_1(M))$ is not
contained in any of the aforementioned closed subgroups, then it is
dense in $\SU$.
\end{rem}
We adopt the following notational conventions:
For a fixed representation $\rho,$ $X\in \pi_1(M),$
we write $X$ for $\rho(X)$
when there is no ambiguity.  A small letter denotes the
trace of the matrix represented by the corresponding
capital letter.

\section{The moduli space of the four-holed sphere}

We first review some results that appear
in \cite{Be1, Go1, Pr2}.
Suppose $M$ is a  three-holed sphere.
Then $\pi_1(M)$ has a presentation:
$$
\langle A, B, C: ABC = I \rangle,
$$
where $A, B,$ and  $C$ represent the homotopy classes of the three boundaries
of $M$.
\begin{prop} \label{prop:Pinon3hole}
\mbox{}

\begin{enumerate}
\item A representation $\rho$ on a three-holed sphere is a
$\Spin$-representation
if and only if $a^2+b^2+c^2-abc-4=0$.
\item A representation $\rho$ on a three-holed sphere is
a $\Pin$-representation and not
a $\Spin$-representation if and only if
$a^2+b^2+c^2-abc-4\neq 0$
and
at least two of the three:
$A,$ $B,$ $AB,$
have zero trace.
\end{enumerate}
\end{prop}
\begin{proof} See \cite{Go1, Pr2}.
\end{proof}

Suppose $M$ is a four-holed sphere.  Then the fundamental
group $\pi_1(M)$ admits a presentation
$$
\langle A, B, C, D : ABCD = I \rangle.
$$
Set $X = AB, Y = BC,$ and  $Z = CA$.
Let $\kappa = (a, b, c, d) \in [-2,2]^4$ be the holonomies on the
boundary.
Then the moduli space $\M_{\kappa}$ is the subspace of $[-2,2]^3$ given
by the equation \cite{Go1, Pr2}
$$
x^2 + y^2 + z^2 + xyz= (ab+cd)x + (ad + bc)y + (ac+bd)z
-(a^2 + b^2 + c^2 + d^2 + abcd -4).
$$
\begin{rem} \cite {Go1}
If two representations in $\M_{\kappa}$
share
$(x,y,z),$ then they are conjugate.
\end{rem}
Let
$$
I_{a,b}= \bigg [\frac { ab - \sqrt{(a^2-4)(b^2-4)}}2, \frac {ab + \sqrt{(a^2-4)
(b^2-4)}}2 \bigg ].
$$
If $I_{a,b} \cap I_{c,d} \neq \emptyset$, then $\M_{\kappa}$ is a
(possibly degenerate) topological sphere (see Figure 1).
%



The mapping class group $\Gamma$ of the 4-holed sphere is generated
by three Dehn twists $\tau_X, \tau_Y, \tau_Z$ \cite{Go1, Pr2}.
In local coordinates, the actions are
$$\left[
\begin{array}{c}
y\\
z\\
\end{array}
\right]
\stackrel{\tau_X}{\longmapsto}
\left[
\begin{array}{c}
ad+bc-x(ac+bd-xy-z)-y\\
ac+bd- xy-z\\
\end{array}
\right],
$$
$$\left[
\begin{array}{c}
z\\
x\\
\end{array}
\right]
\stackrel{\tau_Y}{\longmapsto}
\left[
\begin{array}{c}
bd+ca-y(ba+cd-yz-x)-z\\
ba+cd-yz-x\\
\end{array}
\right],
$$
$$\left[
\begin{array}{c}
x\\
y\\
\end{array}
\right]
\stackrel{\tau_Z}{\longmapsto}
\left[
\begin{array}{c}
cd+ab- z(cb+ad-zx-y)-x\\
cb+ad- zx - y\\
\end{array}
\right].
$$

\

\centerline{\epsfxsize=4in \epsffile{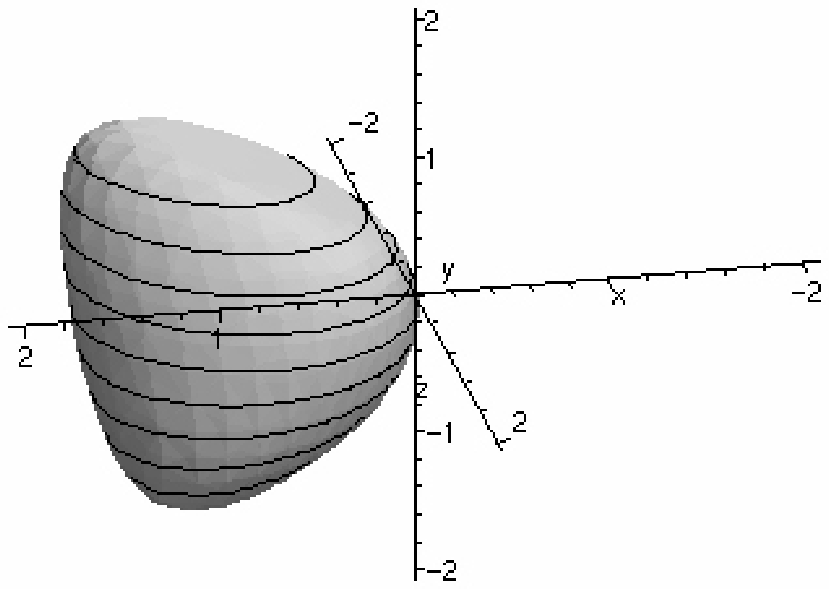}}

\centerline{F{\sc igure} 1: The sphere $\M_\kappa$ for $\kappa =(\sqrt 2, \sqrt 2, \frac 1 2, -\frac 1 2).$}
\

\section{The Pin(2) Representations}

Consider
$$
e^{i \theta} = \begin{pmatrix}
\cos(\theta) & \sin(\theta) \\
-\sin(\theta) & \cos(\theta)
\end{pmatrix}, \ \ \
\iota=\begin{pmatrix}
i & 0 \\
0 & -i
\end{pmatrix}
$$
in $\Pin$.


\begin{prop} \label{prop:spin} Suppose
$\rho \in {\mathcal H}_{(a,b,c,d)}$ with $a,b,c,d \not \in \{\pm 2\}$
and $[\rho] = (x,y,z) \in {\mathcal M}_{\kappa}$.
Then the representation $\rho$
is a $\Spin$-representation if and only if
$x$ is an endpoint of both
$I_{a,b}$ and  $I_{c,d},$
$y$ is an endpoint of both
$I_{b,c}$ and  $I_{a,d},$
and $z$ is an endpoint of both
$I_{a,c}$ and  $I_{b,d}.$
\end{prop}
\begin{proof}
First, suppose that $\rho$ is a $\Spin$-representation. Then, up
to conjugation,
$$
\rho(A) = e^{i \theta_a}, \rho(B) = e^{i \theta_b},
\rho(C) = e^{i \theta_c}, \rho(D) = e^{i \theta_d},
$$
where $\theta_a + \theta_b + \theta_c + \theta_d = 0$.
The  endpoints of $I_{a,b}$ are given by
$$\frac 1 2 (ab \pm \sqrt{(4-a^2)(4-b^2)})$$
$$=\cos(\theta_a+\theta_b)+\cos(\theta_a-\theta_b) \pm \frac 1 2 \sqrt{(4-4\cos^2(\theta_a))(4-4\cos^2(\theta_b))}$$
$$=\cos(\theta_a+\theta_b)+\cos(\theta_a-\theta_b) \pm |2\sin(\theta_a)\sin(\theta_b) |$$
$$=\cos(\theta_a+\theta_b)+\cos(\theta_a-\theta_b)\pm |\cos(\theta_a-\theta_b)-\cos(\theta_a+\theta_b)|.$$
This implies that an endpoint of $I_{a,b}$ is equal to
$2\cos(\theta_a+\theta_b).$
Similarly,
an endpoint of $I_{c,d}$ is equal to $2\cos(\theta_c+\theta_d)$
which is equal to $2\cos(\theta_a+\theta_b)=x.$ Thus $x$ is equal
to an  endpoint of both $I_{a,b}$ and $I_{c,d}.$ A similar
argument shows that $y$ must be an endpoint of $I_{b,c}$ and
$I_{a,d},$ and also  $z$ must be an endpoint of  $I_{a,c}$ and
$I_{b,d}.$

To prove the converse, suppose that
$\rho$ is such that
$x$ is an endpoint of both
$I_{a,b}$ and  $I_{c,d},$
$y$ is an endpoint of both
$I_{b,c}$ and  $I_{a,d},$
and
 $z$ is an endpoint of both
$I_{a,c}$ and  $I_{b,d}.$
Then $2x = ab \pm \sqrt{(4-a^2)(4-b^2)}$ which implies that
$$
4x^2 = a^2b^2 +16-4a^2-4b^2+a^2b^2 \pm 2ab\sqrt{(4-a^2)(4-b^2)}
$$
$$= a^2b^2 +16-4a^2-4b^2+a^2b^2 \pm 2ab(\pm(2x-ab)).$$
Hence
$$x^2 +a^2 +b^2 -xab= 4$$ which implies that $\rho$ is
a $\Spin$-representation on the three-holed sphere $(A,B,X)$ by
Proposition \ref{prop:Pinon3hole}. Similarly, $(C,D,X),$
$(A,C,Z),$ $(B,D,Z),$ $(A,D,Y),$ and $(B,C,Y)$ are all
$\Spin$-representations.
As $A,B,C,$ and $D$ all pairwise commute, we have that $\rho$ is a
$\Spin$-representation on the entire four-holed sphere.
\end{proof}

\begin{prop} \label{prop:pin}
Let $\rho \in {\mathcal H}_{\kappa}$ and
$[\rho] = (x,y,z) \in {\mathcal M}_{\kappa}$.  Suppose $\rho$ is
a $\Pin$-representation but not a $\Spin$-representation
then one of the following two conditions holds:
\begin{enumerate}
\item $\kappa = (0,0,0,0),$
\item $\kappa = (0,0,c,d),$ where $y=0$ and $z=0,$ along with the five other
symmetric cases.
\end{enumerate}
If $\rho$ satisfies one of the two conditions above, then $\rho$ is
a $\Pin$-representation.
\end{prop}

\begin{proof}
Let $\rho$ be a $\Pin$-representation but not a
$\Spin$-representation. Then at least one of $A,B,C,$ or $D$ must
be in $\Spinm.$ However, since $ABCD=I,$ at least two of $A,B,C,$
or $D$ must be in $\Spinm.$ Suppose
 $A,B \in \Spinm.$
If $C\in \Spinm,$ then $D\in \Spinm,$ then we obtain $\kappa = (0,0,0,0).$
If $C\in \Spin,$ then $D\in \Spin,$ which implies that
$AC, BC \in \Spinm,$  i.e., $y=z=0.$

Now consider
$$
A = \iota, B = -\iota e^{i \theta}
$$
which are contained in a $\Pin$ subgroup.

Case 1:
Let $\rho \in {\mathcal H}_{\kappa}$ with $\kappa = (0,0,0,0)$
with $x,y,z$ satisfying the equation
$x^2+y^2+z^2+xyz=4.$
We construct a $\Pin$-representation conjugate to $\rho$ by
setting $x = 2\cos \theta$ (in $A$ and $B$ above)
and setting $C$ equal to one of
$e^{\pm i \psi } \iota$, where
$z= -2\cos \psi.$ As
$CA = -e^{\pm i \psi}$ and
$Y= BC$ is either $e^{i(\theta+\psi)}$ or $e^{i(\theta-\psi)}$
whose traces are the two solutions of
$x^2+y^2+z^2+xyz=4$ for fixed $x$ and $z.$
Therefore, this $\Pin$-representation is
conjugate to $\rho.$

Case 2:
Let $\rho \in {\mathcal H}_{\kappa}$ with $\kappa = (0,0,c,d)$
with $y=z=0.$
Thus $x,c,d$ satisfy:
$x^2 =cdx-c^2-d^2 +4$ implying that
$\rho$ restricted to $(X,C,D)$ is a $\Spin$-representation by
Proposition \ref{prop:Pinon3hole}.
We construct a $\Pin$-representation conjugate to $\rho$ by
setting $x = 2\cos \theta$ (in $A$ and $B$ above)
and setting $C$ to be $e^{i\psi}$
and $D= e^{-i(\psi+\theta)}$.
As the traces of $Y=BC$ and $Z=AC$ are zero,
this $\Pin$-representation is conjugate to $\rho.$
\end{proof}

Propositions \ref{prop:spin} and \ref{prop:pin} provide
a complete characterization of the $\Pin$-representation classes.

\section{Examples}
A direct computation shows that the traces of elements in the groups
${\mathcal C},{\mathcal D}$
are in the set
$$S = \{0, \pm 1, \pm \sqrt 2, \pm \frac{\sqrt{5} + 1}{2},
\pm \frac{\sqrt{5} - 1}{2}, \pm 2 \}.$$
Let $F$ be the set of
$\kappa = (a,a,c,-c) \in [-2,2]^4$ satisfying the following conditions:
\begin{enumerate}
\item $a^2 + c^2 < 4$,
\item $a \neq 0$ and $c \neq 0$,
\item $a \not\in S$ or $c \not\in S$.
\end{enumerate}
Consider the space $\M_\kappa$ with $\kappa \in F$.
A direct computation shows
$$\O = \{(a^2-2,0,0), (2-c^2,0,0)\} \subset \M_\kappa$$
is $\Gamma$-invariant.
By condition 1,
$$
I_{a,a} \cap I_{c,-c} = [a^2-2, 2] \cap [-2,2-c^2] = [a^2-2,2-c^2]
$$
is a closed interval.  Again by condition 1, $a,c \neq \pm 2$.
Hence Proposition \ref{prop:spin} implies that elements in $\O$ do
not correspond to $\Spin$-representations.  By condition 2, $a,c
\neq 0$, so the elements in $\O$ do not correspond to
$\Pin$-representations by Proposition \ref{prop:pin}.  Finally, by
condition 3, they do not correspond to ${\mathcal C}, {\mathcal
D}$-representations. Thus, by Remark \ref{rem:dense}, the elements
in the discrete orbit $\O$ correspond to representations with
dense images in $\SU$.  This proves Theorem \ref{thm:4}.


Figure 1 shows one such case with
$\kappa = (\sqrt 2, \sqrt 2, \frac 1 2, -\frac 1 2)$.  The special
orbit $\O$ consists of the two points that are intersections of the $x$-axis
with $\M_\kappa$, i.e.
$\O=\{ (0,0,0), (\frac 7 4, 0, 0)\}.$
Below is a representation in the conjugacy class
$(0,0,0) \in \O \subset \M_\kappa$:
$$A = B =
\left[
\begin{array}{cc}
\frac {\sqrt 2} 2 + \frac {\sqrt 2} 2 i & 0 \\
0 & \frac {\sqrt 2} 2 - \frac {\sqrt 2} 2 i\\
\end{array}
\right]  \ \ \mbox{and} \ \ \
C = -D =
\left[
\begin{array}{cc}
\frac 1 4 + \frac 1 4 i  &  \frac {\sqrt{14}} 4 \\
 -\frac {\sqrt{14}} 4   &  \frac 1 4 - \frac 1 4 i \\
\end{array}
\right]. \ \ $$

\end{document}